\documentclass{elsart}
\usepackage{amsfonts, amssymb}
\usepackage{amsmath}
\usepackage{latexsym}
\usepackage{verbatim}
\usepackage{url}
\usepackage{graphicx, color}  
\usepackage{enumerate}

\def\t{\mbox{\textbf{\textsf{t}}}}

\newcommand{\NN}{\mathbb{N}}

\def\Alph{\textrm{Alph}}
\def\Ult{\textrm{Ult}}

\newenvironment{remq}{\begin{rem}\rm }{\end{rem}}
\newenvironment{example}{\begin{exmp}\rm }{\end{exmp}}

\def\t1{\mathbf{t}^{(1)}}
\def\bt{\mathbf{t}}
\def\bm{\mathbf{m}}
\def\m1{\mathbf{m}^{(1)}}

\def\bs{\mathbf{s}}
\def\cA{\mathcal{A}}

\def\bu{\mathbf{u}}
\def\bv{\mathbf{v}}
\def\bx{\mathbf{x}}
\def\by{\mathbf{y}}
\def\bz{\mathbf{z}}
\def\rev{\widetilde}
\def\cAstar{\mathcal{A}^*}
\def\cAplus{\mathcal{A}^+}
\def\cAw{\mathcal{A}^\omega}
\def\cAinf{\mathcal{A}^\infty}
\def\cB{\mathcal{B}}
\def\empt{\varepsilon}

\numberwithin{thm}{section}

\journal{European Journal of Combinatorics}

\begin{document}


\begin{frontmatter}

\title{Characterizations of finite and infinite episturmian words via lexicographic orderings}

\author{Amy Glen\corauthref{cor}}
\corauth[cor]{Corresponding author.}
\address{LaCIM, Universit\'e du Qu\'ebec \`a Montr\'eal, C.P. 8888, Succ. Centre-ville,
Montr\'eal, Qu\'ebec, CANADA, H3C 3P8}

\ead{amy.glen@gmail.com}
\author{Jacques Justin}
\address{LIAFA, Universit\'e Paris VII, Case 7014, 2 place Jussieu, 75251 Paris Cedex 05, FRANCE}
\ead{jacjustin@free.fr}

\author{Giuseppe Pirillo}
\address{IASI CNR, Unit\`a di Firenze, Viale Morgagni 67/A, 50134 Firenze, ITALY}
\ead{pirillo@math.unifi.it}

\date{August 26, 2006}

\begin{abstract}
In this paper, we characterize by lexicographic order all {\em finite} Sturmian and episturmian words, i.e., all (finite) factors of such infinite words. Consequently, we obtain a characterization of infinite episturmian words in a {\em wide  sense} (episturmian and {\em episkew} infinite words). That is, we characterize the set of all infinite words whose factors are (finite) episturmian. Similarly, we characterize by lexicographic order all {\em balanced} infinite words over a $2$-letter alphabet; in other words, all Sturmian and {\em skew} infinite words, the factors of which are (finite) Sturmian.
\end{abstract}


\end{frontmatter}


\section{Introduction}

The family of {\em episturmian words} is an interesting natural generalization of the well-known 
{\em Sturmian words} (a particular class of binary infinite words) to an arbitrary finite alphabet, introduced by Droubay, Justin, and Pirillo  \cite{xDjJgP01epis} (also see  \cite{aG05powe,jJgP02epis,jJgP04epis,jJlV00retu} for instance). Episturmian words share many properties with Sturmian words and include the well-known \emph{Arnoux-Rauzy sequences}, the study of which began in  \cite{pAgR91repr} (also see \cite{jJgP02onac,rRlZ00agen} for example). 

In this paper, we characterize by {\em lexicographic order} all {\em finite} Sturmian and episturmian words, i.e., all (finite) factors of such infinite words. Consequently, we obtain a characterization of episturmian words in a {\em wide sense} (episturmian and {\em episkew} infinite words).  That is, we characterize the set of all infinite words whose factors are (finite) episturmian. Similarly, we characterize by lexicographic order all {\em balanced} infinite words over a $2$-letter alphabet; in other words, all Sturmian and {\em skew} infinite words, the factors of which are (finite) Sturmian.

To any infinite word $\bt$ we can associate two infinite words $\min(\bt)$ and $\max(\bt)$ such that any prefix of $\min(\bt)$ (resp.~$\max(\bt)$) is the \emph{lexicographically} smallest (resp.~greatest)  amongst the factors of $\bt$ of the same length (see \cite{gP05ineq} or Section \ref{S:words}). Our main results in this paper extend recent work by Pirillo \cite{gP05ineq,gP05mors}, Justin and Pirillo \cite{jJgP02onac}, and Glen \cite{aG06acha}. In the first of these papers, Pirillo proved that, for infinite words $\bs$ on a 2-letter alphabet $\{a,b\}$ with $a<b$, the inequality $a\bs \leq \min(\bs) \leq \max(\bs) \leq b\bs$ characterizes {\em standard Sturmian words} (both aperiodic and periodic). Similarly, an infinite word $\bs$ on a finite alphabet $\cA$ is {\em standard episturmian} if and only if, for any letter $a \in \cA$ and lexicographic order $<$ satisfying $a = \min(\cA)$, we have 
\begin{equation} \label{eq:gP05ineq}
a\bs \leq \min(\bs). 
\end{equation} 
Moreover, $\bs$ is a {\em strict} standard episturmian word (i.e., a {\em standard} Arnoux-Rauzy sequence \cite{pAgR91repr,rRlZ00agen}) if and only if \eqref{eq:gP05ineq} holds with strict equality \cite{jJgP02onac}. In a similar spirit, Pirillo \cite{gP05mors} very recently defined \emph{fine words} over two letters; that is, an infinite word $\bt$ over a 2-letter alphabet $\{a,b\}$ ($a < b$) is said to be \emph{fine} if  $(\min(\bt), \max(\bt)) = (a\bs, b\bs)$ for some infinite word $\bs$. These words are characterized in \cite{gP05mors} by showing that fine words on $\{a,b\}$ are exactly the {\em aperiodic Sturmian} and {\em skew} infinite words. 

Glen \cite{aG06acha} recently extended Pirillo's definition of fine words to an arbitrary finite alphabet; that is, an infinite word $\bt$ is \emph{fine} if  there exists an infinite word $\bs$ such that $\min(\bt) = a\bs$ for any letter $a \in$ Alph$(\bt)$ and lexicographic order $<$ satisfying $a = \min(\mbox{Alph}(\bt))$. (Here, Alph$(\bt)$ denotes the \emph{alphabet} of $\bt$, i.e., the set of distinct letters occurring in $\bt$.) These generalized fine words are characterized in \cite{aG06acha}; specifically, it is shown that an infinite word $\bt$ is fine if and only if $\bt$ is either a {\em strict}  episturmian word, or a {\em strict episkew word} (i.e., a particular kind of non-recurrent infinite word, all of whose factors are episturmian).  Here, we prove further that an infinite word $\bt$ is episturmian in the {\em wide sense} (episturmian or episkew) if and only if there exists an infinite word $\bu$ such that $a\bu \leq \min(\bt)$ for any letter $a\in \cA$ and lexicographic order $<$ satisfying $a = \min(\cA)$. This result follows easily from our characterization of finite episturmian words in Section 4.

This paper is organized as follows. Section 2 contains all of the necessary terminology and notation concerning words, morphisms, and Sturmian and episturmian words. In Section 3, we give a number of equivalent definitions of {\em episkew words}, and recall the aforementioned characterizations of {\em fine words}. Then, in Section 4, we prove a new characterization of finite episturmian words, from which a new characterization of finite Sturmian words is an easy consequence. Lastly, in Section 5, we obtain characterizations of episturmian words in the wide sense and balanced binary infinite words, which follow from the main results in Sections 3 and 4.

\section{Preliminaries}

\subsection{Words and morphisms} \label{S:words}

Let $\cA$ denote a finite alphabet. A (finite) \emph{word} is an
element of the \emph{free monoid} $\cAstar$ generated by $\cA$, in
the sense of concatenation. The identity $\empt$ of $\cAstar$ is
called the \emph{empty word}, and the \emph{free semigroup},
denoted by $\cAplus$, is defined by $\cAplus :=
\cAstar\setminus\{\empt\}$. An \emph{infinite word} (or simply 
\emph{sequence}) $\bx$ is a sequence indexed by $\NN$ with values in $\cA$, 
i.e., $\bx = x_0x_1x_2\cdots$, where each $x_i \in \cA$. The set of 
all infinite words over $\cA$ is denoted by $\cAw$, and we define $\cAinf := \cAstar \cup \cAw$.

If $w = x_{1}x_{2}\cdots x_{m} \in \cAplus$, each 
$x_{i} \in \cA$, the \emph{length} of $w$ is $|w| = m$ and we denote by 
$|w|_a$ the number of occurrences of a letter $a$ in $w$. (Note that 
$|\empt| = 0$.) The \emph{reversal} $\rev{w}$ of $w$ is given by $\rev{w} = x_{m}x_{m-1}\cdots x_{1}$, and if $w = \rev{w}$, then $w$ is called a \emph{palindrome}.

An infinite word $\bx = x_0x_1x_2\cdots$, each $x_i \in \cA$, is said to be \emph{periodic} (resp.~\emph{ultimately periodic}) with period $p$ if $p$ is the smallest positive integer such that $x_i = x_{i+p}$ for all $i \in \NN$ (resp.~for all $i \geq m$ for some $m\in \NN$). If $u$, $v \in \cA^+$, then $v^\omega$ (resp.~$uv^\omega$) denotes the periodic (resp.~ultimately periodic) infinite word $vvv\cdots$ (resp.~$uvvv\cdots$) having $|v|$ as a period. 

A finite word $w$ is a \emph{factor} of $z \in \cAinf$ if $z = uwv$ for some $u \in \cAstar$, $v \in \cAinf$. Further, $w$ is called a \emph{prefix} (resp.~\emph{suffix}) of $z$ if $u = \empt$ (resp.~$v = \empt$).

  An infinite word $\bx \in \cAw$ is called a \emph{suffix} of $\bz \in \cAw$ if there exists a word $w \in \cA^*$ such that $\bz = w\bx$. A factor $w$ of a word $z \in \cAinf$ is \emph{right} (resp.~\emph{left}) \emph{special} if $wa$, $wb$ (resp.~$aw$, $bw$) are factors of $z$
for some letters $a$, $b \in \cA$, $a \ne b$.

For any word $w \in \cAinf$, $F(w)$ denotes the set of all its
factors, and $F_n(w)$ denotes the set of all factors of
$w$ of length $n \in \NN$, i.e., $F_n(w) := F(w)
\cap \cA^n$ (where $|w|\geq n$ for $w$ finite). Moreover, the \emph{alphabet} of $w$ is Alph$(w)
:= F(w) \cap \cA$ and, if $w$ is infinite, we denote by Ult$(w)$ the set of
all letters occurring infinitely often in $w$.  Two infinite words $\bx$, $\by \in \cAw$ are said to be \emph{equivalent} if $F(\by) = F(\bx)$, i.e., if $\bx$ and $\by$ have the same set of factors. A factor of an infinite word $\bx$ is \emph{recurrent} in $\bx$ if it occurs infinitely many times in $\bx$, and $\bx$ itself is said to be \emph{recurrent} if all of its factors are recurrent in it.

Suppose the alphabet $\cA$ is totally ordered by the relation $<$. Then we 
can totally order $\cA^+$ by the \emph{lexicographic order} $<$, 
defined as follows. Given two words $u$, $v \in \cA^+$, we have $u
< v$ if and only if either $u$ is a proper prefix of $v$ or $u =
xau^\prime$ and $v = xbv^\prime$, for some $x$, $u^\prime$,
$v^\prime \in \cAstar$ and letters $a$, $b$ with $a < b$. This is
the usual alphabetic ordering in a dictionary, and we say that $u$
is \emph{lexicographically less} than $v$. This notion 
naturally extends to $\cAw$, as follows. Let $\bu =
u_0u_1u_2\cdots$ and $\bv = v_0v_1v_2\cdots$, where $u_j$, $v_j
\in \cA$. We define $\bu < \bv$ if there exists an index $i\geq0$
such that $u_j = v_j$ for all $j=0,\ldots, i-1$ and $u_{i} < v_{i}$.
Naturally, $\leq$ will mean $<$ or $=$.

Let $w \in \cAinf$ and let $k$ be a positive integer. We denote by $\min(w | k)$ (resp.~$\max(w | k)$) the lexicographically smallest (resp.~greatest) factor of $w$ of length $k$ for the given order (where $|w|\geq k$ for $w$ finite). If $w$ is infinite, then it is clear that $\min(w | k)$ and $\max(w | k)$ are prefixes of the respective words $\min(w | k+1)$ and $\max(w | k+1)$. So we can define, by taking limits, the following two infinite words (see \cite{gP05ineq})
\[
  \min(w) = \underset{k\rightarrow\infty} {\lim}\min(w | k) \quad \mbox{and} \quad 
  \max(w) = \underset{k\rightarrow\infty}{\lim}\max(w | k).
\]  

The \emph{inverse} of $w \in \cAstar$, written $w^{-1}$, is
defined by $ww^{-1} = w^{-1}w = \empt$. It must be emphasized that
this is merely formal notation, i.e., for $u, v, w \in \cAstar$, the
words $u^{-1}w$ and $wv^{-1}$ are defined only if $u$ (resp.~$v$)
is a prefix (resp.~suffix) of $w$.

A \emph{morphism on} $\cA$ is a map $\psi: \cAstar \rightarrow
\cAstar$ such that $\psi(uv) = \psi(u)\psi(v)$ for all $u, v \in
\cAstar$.  It is uniquely determined by its image on the alphabet
$\cA$. The action of morphisms on $\cA^*$ naturally extends to infinite words; that is, if 
$\bx = x_0x_1x_2 \cdots \in \cAw$, then $\psi(\bx) = \psi(x_0)\psi(x_1)\psi(x_2)\cdots$. 

In what follows, we shall assume that $\cA$ contains two or more letters.

\subsection{Sturmian words} \label{SS:Sturmian}

Sturmian words admit several equivalent definitions and have numerous characterizations; for instance, they can be characterized by their palindrome or return word structure
\cite{xDgP99pali,jJlV00retu}. A particularly useful definition of Sturmian words is the following.

\newpage
\begin{defn}
An infinite word $\bs$ over $\{a,b\}$ is \textbf{\em Sturmian} if there exist real numbers $\alpha$, $\rho \in [0,1]$  such that $\bs$ is equal to one of the following two infinite words:
\[
  s_{\alpha,\rho}, ~s_{\alpha,\rho}^{\prime}: \NN \rightarrow \{a,b\} 
\]
defined by
\[
 \begin{matrix}
  &s_{\alpha,\rho}(n) = \begin{cases}
                        a    &\mbox{if} ~\lfloor(n+1)\alpha + \rho\rfloor -
                        \lfloor n\alpha + \rho\rfloor = 0, \\
                        b    &\mbox{otherwise};
                       \end{cases} \\
  &\qquad \\
  &s_{\alpha,\rho}^\prime(n) = \begin{cases}
                        a    &\mbox{if} ~\lceil(n+1)\alpha + \rho\rceil -
                        \lceil n\alpha + \rho\rceil = 0, \\
                        b    &\mbox{otherwise}.
                       \end{cases}
 \end{matrix} \qquad (n \geq 0)
\]
Moreover, $\bs$ is said to be  \textbf{\em standard Sturmian} if $\rho = \alpha$. 
\end{defn}

\begin{remq}  A Sturmian word of {\em slope} $\alpha$ is:
\begin{itemize}
\item \emph{aperiodic} (i.e., not ultimately periodic) if $\alpha$ is irrational;
\item \emph{periodic} if $\alpha$ is rational.
\end{itemize}

Nowadays, for most authors, only the aperiodic Sturmian words are considered to be `Sturmian'. In several of our previous papers (see \cite{aG06acha,jJgP97deci,jJgP04epis,gP99anew,gP05mors} for instance), we have referred to aperiodic Sturmian words as `proper Sturmian' to highlight the fact that such Sturmian words correspond to the most common sense of `Sturmian' now. In the present paper, the term `Sturmian' will refer to both aperiodic and periodic Sturmian words.
\end{remq}

\begin{defn} \label{D:balance} 
A finite or infinite word $w$ over $\{a,b\}$ is said to be \textbf{\em balanced} if, for any factors $u$, $v$ of $w$ with $|u| = |v|$, we have $||u|_{b} - |v|_{b}| \leq 1$ (or equivalently $||u|_{a} - |v|_{a}| \leq 1$). 
\end{defn}

In the pioneering paper \cite{gHmM40symb}, balanced infinite words over a $2$-letter alphabet  are called `Sturmian trajectories' and belong to three classes: aperiodic Sturmian; periodic Sturmian; and non-recurrent infinite words that are ultimately periodic (but not periodic), called {\em skew words}. That is, the family of balanced infinite words consists of the (recurrent) Sturmian words and the (non-recurrent) skew infinite words, all of whose factors are balanced.

It is important to note that a finite word is \emph{finite Sturmian} (i.e., a factor of some Sturmian word) if and only if it is balanced \cite{jBpS02stur}. Accordingly, balanced infinite words are precisely the infinite words whose factors are finite Sturmian. In Section \ref{S:infinite}, we will generalize this concept by showing that the set of all infinite words whose factors are {\em finite episturmian} consists of the (recurrent) episturmian words and the (non-recurrent) {\em episkew} infinite words (see Propositions \ref{P:skew} and \ref{P:wide}, to follow). 

For a comprehensive introduction to Sturmian words, see for instance \cite{jAjS03auto,jBpS02stur,nP02subs} and references therein. Also see \cite{aHrT00char,gP05mors} for further work on skew words.

\subsection{Episturmian words} \label{S:episturmian}

For episturmian words and morphisms\footnotemark[1] we use the same terminology and notation as in \cite{xDjJgP01epis,jJgP02epis,jJgP04epis}. \footnotetext[1]{In \cite{jJgP02epis}, Section 5.1 is incorrect and should be ignored.}

An infinite word $\bt \in \cAw$ is \emph{episturmian} if $F(\bt)$ is closed under
reversal and $\bt$ has at most one right (or equivalently left) special factor of each length. Moreover, an episturmian word is \emph{standard} if all of its left special factors are prefixes of
it. Sturmian words are exactly the episturmian words over a 2-letter alphabet. 

\noindent\textit{Note.} Episturmian words are recurrent \cite{xDjJgP01epis}.

Standard episturmian words are characterized in \cite{xDjJgP01epis} using the concept of the \emph{palindromic right-closure} $w^{(+)}$ of a finite word $w$, which is the (unique) shortest palindrome having $w$ as a prefix (see \cite{aD97stur}).   
Specifically, an infinite word $\bt \in \cAw$ is standard episturmian if and only if  there exists an infinite word $\Delta(\bt) = x_1x_2x_3\ldots$, each $x_i
\in \cA$, called the \emph{directive word} of $\bt$, such that the infinite sequence of palindromic prefixes $u_1 =
\empt$, $u_2$, $u_3$, $\ldots$ of $\bt$ (which exists by results in \cite{xDjJgP01epis}) is given by 
\begin{equation} \label{eq:02.09.04}
  u_{n+1} = (u_nx_n)^{(+)}, \quad 
  n \in \NN^+.
\end{equation}
\textit{Note.} An equivalent way of constructing the sequence $(u_n)_{n\geq 1}$ is via the `hat operation' \cite[Section III]{rRlZ00agen}.

Let $a \in \cA$ and denote by $\psi_a$ the morphism on $\cA$ defined by 
\[
\psi_a : \left\{\begin{array}{lll}
                a &\mapsto &a \\
                x &\mapsto &ax \quad \mbox{for all $x \in \cA\setminus\{a\}$}.
                \end{array} \right.
\]
Together with the permutations of the alphabet, all of the morphisms $\psi_a$ generate by composition the monoid of {\em epistandard morphisms} (`epistandard' is an elegant shortcut for `standard episturmian' due to Richomme \cite{gR03conj}). The submonoid generated by the $\psi_a$ only is the monoid of {\em pure epistandard morphisms}, which includes the {\em identity morphism} Id$_{\cA} =$ Id, and consists of all the \emph{pure standard (Sturmian)  morphisms} when $|\cA|=2$. 

\begin{remq} \label{R:separating} If $\bx = \psi_a(\by)$ or $\bx = a^{-1}\psi_a(\by)$ for some $\by \in \cAw$ and $a \in \cA$, then the letter $a$ is said to be \emph{separating for $\bx$} and its factors; that is, any factor of $\bx$ of length $2$ contains the letter $a$.
\end{remq}

Another useful characterization of standard episturmian words is the following (see \cite{jJgP02epis}). An infinite word $\bt \in \cAw$ is standard episturmian with directive word $\Delta(\bt) = x_1x_2x_3\cdots$ ($x_i \in \cA$) if and only if there exists an infinite sequence of infinite words
$\bt^{(0)} = \bt$, $\bt^{(1)}$, $\bt^{(2)}$, $\ldots$ such that $\bt^{(i-1)} = \psi_{x_i}(\bt^{(i)})$ for all $i \in \NN^+$. Moreover, each $\bt^{(i)}$ is a standard episturmian word with directive word 
$\Delta(\bt^{(i)}) = x_{i+1}x_{i+2}x_{i+3}\cdots$, the \emph{$i$-th shift} of $\Delta(\bt)$. 

To the prefixes of the directive word $\Delta(\bt) = x_1x_2\cdots$, we associate the morphisms 
\[
  \mu_0 := \mbox{Id}, \quad \mu_n := \psi_{x_1}\psi_{x_2}\cdots\psi_{x_n}, \quad n \in \NN^+, 
\]
and define the words 
\[
  h_n := \mu_n(x_{n+1}), \quad n \in \NN, 
\]
which are clearly prefixes of $\bt$.  For the palindromic prefixes $(u_i)_{i\geq 1}$ given by \eqref{eq:02.09.04}, we have the following useful formula \cite{jJgP02epis}
\[
  u_{n+1} = h_{n-1}u_{n};
\]
whence, for $n > 1$ and $0 < p < n$, 
\begin{equation} \label{eq:u_n}
  u_n = h_{n-2}h_{n-3}\cdots h_1h_0 = h_{n-2}h_{n-3}\cdots h_{p-1}u_p.
\end{equation}

\noindent\textit{Note.}  
Evidently, if a standard episturmian word $\bt$ begins with the letter $x \in \cA$, then $x$ is separating for $\bt$ (see \cite[Lemma 4]{xDjJgP01epis}).

\subsubsection{Strict episturmian words}

A standard episturmian word $\bt\in \cA^\omega$, or any equivalent (episturmian) word,   
is said to be \emph{$\cB$-strict} (or $k$-\emph{strict} if $|\cB|=k$, or {\em strict} if $\cB$ is understood) if 
Alph$(\Delta(\bt)) =$ Ult$(\Delta(\bt)) = \cB \subseteq \cA$.  In particular, a standard episturmian word 
over $\cA$ is $\cA$-strict if every letter in $\cA$ occurs infinitely many times in its directive word. 
The $k$-strict episturmian words have complexity $(k - 1)n+1$ for 
each $n \in \NN$; such words are exactly the $k$-letter Arnoux-Rauzy  
sequences. In particular, the $2$-strict episturmian words correspond to the aperiodic Sturmian words. The strict standard episturmian words are precisely the standard (or characteristic) Arnoux-Rauzy sequences.

\section{Episkew words} \label{S:skew}

Recall that a finite word $w$ is said to be {\em finite Sturmian} (resp.~\emph{finite episturmian}) if $w$ is a factor of some infinite Sturmian (resp.~episturmian) word. When considering factors of infinite episturmian words, it suffices to consider only the strict standard ones (i.e., characteristic Arnoux-Rauzy sequences). Indeed, for any factor $u$ of an episturmian word, there exists a strict standard episturmian word also having $u$ as a factor.  
Thus, finite episturmian words are exactly the {\em finite Arnoux-Rauzy words} considered by Mignosi and Zamboni \cite{fMlZ02onth}.

In this section, we define {\em episkew words}, which were alluded to (but not explicated) in the recent paper \cite{aG06acha}. The following proposition gives a number of equivalent definitions of such  infinite words.

\noindent\textbf{Notation:} Denote by $\bv_p$ the prefix of length $p$ of a given infinite word $\bv$.

\begin{prop} \label{P:skew} An infinite word $\bt$ with {\em Alph}$(\bt) = \cA$ is \textbf{\em episkew} if equivalently:
\begin{enumerate}[{\em (i)}]
\item $\bt$ is non-recurrent and all of its factors are (finite) episturmian;

\item there exists an infinite sequence $(\bt^{(i)})_{i\geq0}$ of non-recurrent infinite words and a directive word $x_1x_2x_3\cdots$ $(x_i \in \cA)$ such that $\bt^{(0)} = \bt$, $\ldots$~, $\bt'^{(i-1)} = \psi_{x_i}(\bt^{(i)})$, where  $ \bt'^{(i-1)}= \bt^{(i-1)}$ if $ \bt^{(i-1)}$ begins with $x_i$ and ${\bt^\prime}^{(i-1)} = x_i\bt^{(i-1)}$ otherwise;
\item  there exists a letter $x \in \cA$ and a standard episturmian word $\bs$ on $\cA\setminus\{x\}$ such that $\bt = v \mu(\bs)$, where  $\mu$ is a pure epistandard morphism on $\cA$ and $v$ is a non-empty suffix of $\mu(\rev{\bs_p}x)$ for some $p \in \NN$.
\end{enumerate}
Moreover, $\bt$ is said to be a \textbf{\em strict episkew word} if $\bs$ is strict on $\cA\setminus\{x\}$, i.e., if each letter in $\cA\setminus\{x\}$ occurs infinitely often in the directive word $x_1x_2x_3\cdots$.
\end{prop}
\vspace{-15pt}\begin{pf}
\noindent (i) $\Rightarrow$ (ii): Since all of the factors of $\bt$ are finite episturmian, there exists a letter,  $x_1$ say, that is separating for $\bt$. If $\bt$ does not begin with $x_1$, consider $\bt' = x_1\bt$; otherwise consider $\bt' = \bt$. Then, $\bt' = \psi_{x_1}(\bt^{(1)})$ for some $\bt^{(1)} \in \cAw$. Continuing in this way, we obtain infinite words $\bt^{(2)}$, $\bt'^{(2)}$, $\bt'^{(3)}$, $\bt^{(3)}$, $\ldots$  with $\bt'^{(i-1)}$ as in the statement. Clearly, if some $\bt^{(i)}$ is recurrent then $\bt$ is also recurrent, in which case $\bt$ is episturmian by \cite[Theorem 3.10]{jJgP02epis}. Thus all of the $\bt^{(i)}$ are non-recurrent. 

\medskip
 
\noindent (ii) $\Rightarrow$ (iii): We proceed by induction on $| \cA|$. The starting point of the induction  (i.e., $|\cA|=2$) will be considered later. 

Let $\Delta := x_1x_2x_3 \cdots$. If $\cA= \Ult(\Delta)$ then any letter in $\cA$ is separating for infinitely many  $\bt^{(i)}$, thus is recurrent in all $\bt^{(i)}$. Consider any factor $w$ of $\bt$. As $ |\Ult(\Delta )| >1$, we easily see that $w$ is a factor of $\psi_{x_1}\psi_{x_2}\cdots\psi_{x_q}(x)$ for some $q$ and letter $x$. Hence $w$ is recurrent in $\bt$ and it follows that $\bt$ itself is recurrent; a contradiction. Thus, there exists a letter $x$ in $\cA$ and some minimal $n$ such that $x$ is not recurrent in $\bt^{(n)}$. Two cases are possible: \smallskip

\noindent {\em Case $1$}: $x$ does not occur in $\bt^{(n)}$. Then $ |\Alph(\bt^{(n)})| < |\cA|$; whence, by induction, $\bt^{(n)}$ has the desired form and clearly $\bt$ also has the desired form. More precisely, if we let $\cB := \cA\setminus\{x\}$, then $\bt^{(n)} = \hat v \lambda(\bs)$ where $\bs$ is a standard episturmian word on $\cB\setminus\{y\}$ for some letter $y \ne x$, $\lambda$ is a pure epistandard morphism on $\cB$, and $\hat v$ is a non-empty suffix of $\lambda(\rev{ \bs _q}y)$ for some $q \in \NN$. It easily follows that $\bt = v\mu(\bs)$ where $\bs$ is a standard episturmian word on $\cA\setminus\{y\}$, $\mu$ is a pure epistandard morphism on $\cA$, and $v$ is a non-empty suffix of $\mu(\rev{\bs_p}y)$ for some $p \in \NN$. \smallskip

\noindent {\em Case $2$}: $x$ occurs in $\bt^{(n)}$. We now show that $x$ occurs exactly once in $\bt^{(n)}$. 

Suppose on the contrary that $x$ occurs at least twice in $\bt^{(n)}$. Then,    
since $x_{n+1}$ is separating for $\bt^{(n)}$, we have $x w^{(n)}x \in F(\bt)$ for some non-empty word $w^{(n)}$ for which $x_{n+1}$ is separating, and the first and last letter of $w^{(n)}$ is $x_{n+1}$ (that is,
 $w^{(n)}x = \psi_{x_{n+1}}(w^{(n+1)}x)$,
 where $w^{(n+1)} = \psi_{x_{n+1}}^{-1}(w^{(n)}x_{n+1}^{-1})$).
 Using the fact that $ |w^{(n)}x|= 2 |w^{(n+1)}x| - |w^{(n+1)}x|_{x_{n+1}}$, we see that $|w^{(n+1)}| < |w^{(n)}|$. So, continuing the above procedure, we obtain infinite words $\bt^{(n+1)}$, $\bt^{(n+2)}$, $\ldots$~  containing similar shorter factors $xw^{(n+1)}x$, $x w^{(n+2)}x$, $\ldots$ until we reach $\bt^{(q)}$, which contains $xx$. But this is impossible because the letter $x_{q+1} \ne x$ is separating for $\bt^{(q)}$.  Therefore $\bt^{(n)}$ contains only one occurrence of $x$ and we have 
$$\bt^{(n)} = ux\bs^{(n)} \quad \mbox{for some $u \in (\cA\setminus\{x\})^*$ and $\bs^{(n)} \in (\cA\setminus\{x\})^\omega$}.$$ 
Now, as $x$ is never separating for $\bt^{(j)}$, $j \geq n$, we can write  $\bt^{(n+j)} = u^{(j)}x\bs^{(n+j)}$  for some $u^{(j)}$, $\bs^{(n+j)}$, and we have $ \bs^{(n+j-1)} = \psi_{x_{n+j}}(\bs^{(n+j)})$, $j>0$. It follows by the Preliminaries (Section \ref{S:episturmian}) that $\bs^{(n)}$ is a (recurrent) standard episturmian word.
 
Now we study the factor $u$ preceding $x$ in $\bt^{(n)}$. Let $u^\prime = x_{n+1}u$ if $u$ does not begin with $x_{n+1}$; otherwise let $u^\prime = u$. Then $u^\prime x$ is a prefix of ${\bt'}^{(n)}$. Moreover, since $x_{n+1}$ is separating for $u^\prime x$,  we have $u^\prime x = \psi_{x_{n+1}}(u^{(1)}x)$ where $u^{(1)} = \psi_{x_{n+1}}^{-1}(u^\prime x_{n+1}^{-1})$. Hence $\bt^{(n+1)} = u^{(1)}x \bs^{(n+1)}$, where $x_{n+2}$ is separating for $u^{(1)}x $ (if $u^{(1)} \ne \empt$).  Continuing in this way,  we arrive at the infinite word $\bt^{(q)} = x \bs^{(q)}$ for some $q \geq n$, where $\bs^{(q)}$ is a standard episturmian word on $\cA\setminus\{x \}$.   

Reversing the procedure, we find that 
\[
  \bt^{(n)} = w\bs^{(n)} \quad \mbox{where $w = ux$ is a non-empty suffix of $\psi_{x_{n+1}}\cdots \psi_{x_q}(x)$}.
\]
Suppose $(u_i)_{i\geq1}$ is the sequence of palindromic prefixes of 
\[
\bs = \psi_{x_1}\cdots \psi_{x_n} (\bs^{(n)}) = \mu_n(\bs^{(n)}),
\] 
and the words $(h_i)_{i\geq0}$ are the prefixes  $(\mu_{i}(x_{i+1}))_{i\geq 0}$ of $\bs$. Then, letting $u_i^{(n)}$,  $h_i^{(n)}$, and $\mu_i^{(n)}$ denote the analogous elements for $\bs^{(n)}$, we have  
$$\mu_0^{(n)} = \mbox{Id}, \quad \mu^{(n)}_{i} = \psi_{x_{n+1}}\psi_{x_{n+2}}\cdots \psi_{x_{n+i}} = \mu_n^{-1}\mu_{n+i}$$ and 
\[  
h_0^{(n)} = x_{n+1}, \quad h_i^{(n)} = \mu_{i}^{(n)}(x_{n+1+i}) \quad \mbox{for $i=1$, $2$, $\ldots$~.}
\]
Now, if $u\ne \empt$, then $q \geq 1$, and we have 
\begin{align*}
 \psi_{x_{n+1}}\cdots \psi_{x_q}(x) = \mu_{q-n}^{(n)}(x)  
                        &= \mu_{q-n-1}^{(n)}\psi_q(x_n) \notag\\
                        &= \mu_{q-n-1}^{(n)}(x_qx_n) \notag\\
                        &= h_{q-n-1}^{(n)}\mu_{q-n-1}^{(n)}(x) \notag\\
                        &\qquad \vdots \notag\\
                        &=h_{q-n-1}^{(n)}\cdots h_{1}^{(n)}\mu_0^{(n)}(x_{n+1}x) \notag \\
                        &=h_{q-n-1}^{(n)}\cdots h_1^{(n)}h_0^{(n)}x  
                        =u_{q-n+1}^{(n)}x \quad \mbox{(by \eqref{eq:u_n}).} 
\end{align*}  
Therefore, $w = ux$ where $u$ is a (possibly empty) suffix of the palindromic prefix $u_{q-n+1}^{(n)}$  of $\bs^{(n)}$. That is, $u$ is the reversal of some prefix of $\bs^{(n)}$; in particular 
$$u = \rev{\bs}_p^{(n)} \quad \mbox{for some $p \in \NN$},$$
and hence
\[
 \bt^{(n)} = \rev{\bs}_p^{(n)}x\bs^{(n)}.
\]
So, passing back from $\bt^{(n)}$ to $\bt$, we find that  
\[
  \bt = v\mu_n(\bs^{(n)}) = v\bs \quad \mbox{where $v$ is a non-empty suffix of $\mu_n(\rev{\bs}_p^{(n)}x)$. } \]

It remains to treat the case $| \cA| = 2$. Reasoning as previously we see that for some $n$, $\bt^{(n)} = y^p x y ^\omega$ where $x \neq y \in \cA$; whence the desired form for $\bt$.
  \medskip

\noindent (iii) $\Rightarrow$ (i): It suffices to show that the factors of $\rev{\bs_p}x\bs$ are (finite)  episturmian. This is trivial for factors not containing the letter $x$. Suppose $w$ is a factor containing $x$. Then $w$ is a factor of $u_r x u_r$ where $u_r$ is a long enough palindromic prefix of $\bs$. Thus it remains to show that $u_r x u_r$ is episturmian and this is true because it is $(u_rx)^{(+)}$, which is a palindromic prefix of some standard episturmian word. 
\qed\end{pf} 

\vspace{-10pt}
\begin{remq}
Episkew words on a $2$-letter alphabet are precisely the skew words, defined in Section \ref{SS:Sturmian}.
\end{remq}

\subsection{Fine words}

\begin{defn} An {\bf acceptable pair} is a pair $(a, <)$ where $a$ is a letter and $<$ is a lexicographic order on $\cA^+$ such that $a = \min(\cA).$
\end{defn}

\begin{defn} \cite{aG06acha}
An infinite word $\bt$ on $\cA$ is said to be {\bf fine} if  there exists an infinite word $\bs$ such that $\min(\bt) = a\bs$ for any acceptable pair $(a,<)$.
\end{defn}

\noindent\textit{Note.} 
Since there are only two lexicographic orders on words over a 2-letter alphabet, a fine word $\bt$ over $\{a,b\}$ ($a< b$) satisfies $(\min(\bt), \max(\bt)) = (a\bs, b\bs)$ for some infinite word $\bs$.

Pirillo \cite{gP05mors} characterized fine words over a 2-letter alphabet. Specifically:

\begin{prop} \label{P:gP05mors} Let $\bt$ be an infinite word over $\{a,b\}$. The following properties are equivalent:  
\begin{enumerate}[{\em (i)}]
\item $\bt$ is fine,
\item either $\bt$ is aperiodic Sturmian, or $\bt = v\mu(x)^\omega$ where $\mu$ is a pure standard Sturmian morphism on $\{a,b\}$, and $v$ is a non-empty suffix of $\mu(x^py)$ for some $p \in \NN$ and $x$, $y \in \{a,b\}$ $(x\ne y)$.  \qed
\end{enumerate} 
\end{prop}

In other words, a fine word over two letters is either an aperiodic Sturmian word or an ultimately periodic (but not periodic) infinite word, all of whose factors are Sturmian, i.e., a \emph{skew word} (see Section \ref{SS:Sturmian}). Recently, Glen \cite{aG06acha} generalized this result to infinite words over two or more letters; that is, an infinite word $\bt$ is fine if and only if $\bt$ is either a strict episturmian word or a strict episkew word.

\section{A characterization of finite episturmian words}

Let $w \in \cAinf$ and let $k$ be a positive integer. Recall that $\min(w | k)$ (resp.~$\max(w | k)$) denotes the lexicographically smallest (resp.~greatest) factor of $w$ of length $k$ for the given order (where $|w|\geq k$ for $w$ finite).

\begin{defn} For a finite word $w \in \cA^+$ and a given order, $\min(w)$ will denote $\min(w | k)$ where $k$ is maximal such that all $\min (w | j)$, $j= 1,2, \dots, k$, are prefixes of $\min (w | k)$. In the case $\cA= \{a,b\}$, $\max(w)$ is defined similarly.
\end{defn} 

\begin{example} \label{Ex:min} Suppose $w = baabacababac$. Then, for the orders $b<a<c$ and $b<c<a$ on the $3$-letter alphabet $\{a,b,c\}$: \vspace{-10pt}
\begin{eqnarray*}
\min(w|1) &=& b \\
\min(w|2) &=& ba \\
\min(w|3) &=& bab \\
\min(w|4) &=& baba \\
\min(w|5) &=& babac ~=~ \min(w)
\end{eqnarray*}
\end{example} \vspace{-5pt}
Notice that, in the above example, $\min(w)$ is a suffix of $w$; in fact, this interesting property is true in general, as shown below.

\begin{prop} \label{prop:suffix}
For any finite word $w$ and a given order, $\min(w)$ is a suffix of $w$. Moreover, $\min(w)$ is unioccurrent (i.e., has only one occurrence) in $w$.
\end{prop}
\vspace{-15pt}\begin{pf} If  $ \min (w)$ ($= \min (w | k)$, say) has an occurrence in $w$ that is not a suffix of $w$, then $\min(w | k+1) = \min (w | k)x$ for some letter $x$, contradicting the maximality of $k$. Hence $\min(w)$ occurs just once in $w$ as a suffix.
\qed\end{pf}

\vspace{-10pt}
\noindent\textbf{Notation:}  From now on, it will be convenient to denote by $v_p$ the prefix of length $p$ of a given finite {\em or} infinite word $v$ (where $|v| \geq p$ for $v$ finite). 

In this section, we shall prove the following characterization of finite episturmian words.

 \begin{thm}\label{t3} A finite word $w$ on $\cA$ is episturmian if and only if there exists a finite word $u$ such that, for any acceptable pair  $(a, <)$, we have 
\begin{equation}au_{|m|-1} \le m \label{e2} \end{equation} where $m= \min(w)$ for the considered order.
\end{thm}

The following two lemmas are needed for the proof of Theorem \ref{t3}.

\begin{lem} \label{lem:separating}
If $w$ and $u$ satisfy inequality $\eqref{e2}$ for all acceptable pairs $(a,<)$ and $|\mbox{{\em Alph}}(w)| >1$, then $u$ is non-empty and its first letter is separating for $w$.  
\end{lem}
\vspace{-15pt}\begin{pf} 
Let $a \neq b \in \Alph(w)$ and let $(a,<)$, $(b,<')$ be  two acceptable pairs. As the corresponding two $\min(w)$'s are suffixes of  $w$ (by Proposition \ref{prop:suffix}), they have different lengths; whence $|u| > 0$.

Now we show that the first letter $u_1$ of $u$ is separating for $w$. Indeed, if this is not true, then there exist letters $z$, $z' \in \cA\setminus\{u_1\}$ (possibly equal) such that $zz' \in F(w)$. But $\min(\cA) = z \leq z' < u_1$ for some acceptable pair $(z,<)$, in which case $zz' < zu_1$, contradicting the fact that $zu_1 \leq m_2$.
\qed\end{pf}

\begin{lem} \label{lem:psi_z}
Consider $w$, $w'$ $ \in \cA^*$ and some letter $z \in \cA$. For any given order $<$ on $\cA$:
\begin{enumerate}[{\em (i)}]
\item if $w$ does not end with $z$ and $w= \psi_z(w')$, then 

\[
  \min(w) = \begin{cases}
                    \psi_z(\min(w')) &\mbox{if $\min(w)$ begins with $z$}, \\
                    z^{-1}\psi_z(\min(w')) &\mbox{otherwise};
                    \end{cases}
\]                      
\item if $w$ ends with $z$ and $w= \psi_z(w')z$, then  

\[
  \min(w) = \begin{cases}
                    \psi_z(\min(w'))z &\mbox{if $\min(w)$ begins with $z$}, \\
                    z^{-1}\psi_z(\min(w'))z &\mbox{otherwise}.
                    \end{cases}
\]                      
\end{enumerate}                       
\end{lem}
\vspace{-15pt}\begin{pf}
We denote by $m$, $m'$ the respective words $\min (w)$, $\min (w')$. 

Consider first the simplest case: $w$ does not end with $z$, $m$ begins with $z$. Thus $w= \psi_z (w')$ for some word $w'$ that does not end with $z$. Write $e:= \psi_z (m')$. We have to show that $e=m$.  
Let $k$ be maximal such that $e_i=\min(w | i)$ for $i=1, \dots, k$. Suppose $k<|e|$. Then there exist $x,y \in \cA$, $x>y$, such that $e_{k+1} = e_k x$ and $e_k y \in 
F(w)$. Thus, as $z$ is separating for $w$, $e_k = e_{k-1}z$, with $e_{k-1}=\psi_z (m'_q)$ for some $q$. Since $m$ begins with $z$, $\min(\Alph(w')) = z$ and we have $e_{k+1} = \psi_z(m'_q x)= \psi_z (m'_{q+1})$. Also, if $y \neq z$ then $e_k y= \psi_z (m'_q y)$ with $m'_q y \in F(w')$. If $y=z$, then as $w$ does not end with $z$, $e_k yd = e_{k-1} zyd $ is a factor of $w$ for some letter $d$;  whence again  $m'_q y \in F(w')$. As $x>y$, this contradicts $m'_{q+1}= \min(w' | q+1)$.

Thus $k=|e|$. It suffices now to show that no $ex$, $x \in \cA$, occurs in $w$. Otherwise $ex \in F(w)$. As $m'$ does not end with $z$, also $e$ does not end with $z$, thus $x=z$. So, as $w$ does not end with $z$, $ezy=exy$ occurs in $w$ for some letter $y$, whence  $\psi_z (m' y) \in F(w')$ contradicting the unioccurrence of $m'$ in $w'$.

Now we pass to the most complicated case: $w$ ends with $z$, $m$ does not begin with $z$, $w= \psi_z(w')z$. Letting $e:= z^{-1} \psi_z(m')z$, we need to show that $e=m$. Let $k$ be maximal such that $e_i=\min(w | i)$ for $i=1, \dots, k$. Suppose $k<|e|$. Then, there exist $f \in F(w)$ and $x,y \in \cA$ with $x>y$, such that $e_{k+1} = e_k x$ and $f=e_k y$. As $w$ begins with $z$, clearly $z e_{k+1}, zf \in F(w)$. Also $e_k$ ends with $z$, hence  
$ze_{k+1} = \psi_z (m'_q)zx $ and $zf=\psi_z(m'_q)zy$ for some $q <|m'|$. We distinguish three cases:  $x,y \neq z$; $x=z$; $y=z$. 

The first case leads to $ze_{k+1} = \psi_z (m'_qx)$ and $zf=\psi_z(m'_qy)$;  whence $m'_{q+1} >m'_q y$, contradicting the definition of $m'$. For the case $x=z, y \neq z$, let $m' = m'_q u$, $u \in \cAstar$, and recall that $ze= \psi_z(m')z$. We get $ze=  \psi_z (m'_q) \psi_z (u)z$, thus $\psi_z (u)z$ begins with $zz$, and so  $u$ begins with $z$. Hence $m'_{q+1}= m'_qz=m'_qx$, leading to a contradiction as above. The third case is similar.

Thus $k= |e|$ and it remains to show that no $ex$, $x$ a letter, occurs in $w$. Consider for instance the case $x=z$. Indeed $ez \in F(w)$ implies $z^{-1}\psi_z (m')zz \in F(\psi_z (w')z)$, so 
$\psi_z (m')z \in F(\psi_z (w'))$, whence $m'd \in F(w')$ for some letter $d$; a contradiction.

The other two cases in the lemma have similar proofs.
\qed\end{pf}

\vspace{-5pt}
\begin{example} Let us illustrate the most complicated case when $w$ ends with $z$ and $m$ does not begin with $z$. Let $w' = aa$, $z=b$, $w = babab = \psi_b (w')b $. Then $m'=aa$ and $m= abab= b^{-1}\psi_b (m')b$. 
\end{example}

\vspace{-15pt}\begin{pf*}{Proof of Theorem \ref{t3}}  ONLY IF part: $w$ is finite episturmian, so is a factor of some standard episturmian word $\bs$. By \cite[Proposition 3.2]{gP05ineq} or \cite[Theorem 0.1]{jJgP02onac}, $a\bs\leq \min(\bs)$ for any acceptable pair $(a,<)$. Thus, $m = \min(w)$ trivially satisfies 
\[
  a\bs_{|m|-1} \leq m;
\] 
that is, with $r$ large enough and $u=\bs_r$, inequality \eqref{e2} is satisfied for any acceptable pair $(a,<)$, as required. \smallskip

\noindent IF part: Remark first that if \eqref{e2} is satisfied for some $u$ then it also holds for any $uv$, $v \in \cAstar$. Also, if $a \not \in \Alph (w)$ then \eqref{e2} is trivially satisfied, allowing us to limit our  attention to acceptable pairs $(a,<)$ with $a \in \Alph (w)$.

Let $x:= u_1$, the first letter of $u$. The proof will proceed by induction on $\ell= |w|$. If $w$ is a letter, then $w$ is clearly finite episturmian, i.e., the initial case $|w| = 1$ is trivially true. 

We now distinguish two cases according to whether or not $w$ begins with $x$. 

\noindent \emph{Case $1$}: $w$ begins with $x$. Suppose for instance $w$ does not end with $x$ (the other case is similar). Then, by Lemma \ref{lem:separating},  $w = \psi_x (w')$ for some word $w'$ that does not end with $x$. Further, it follows from Lemma \ref{lem:psi_z} that, for any acceptable pair $(a,<)$, $\min(w) = \psi_x(\min(w'))$ if $x=a$ (resp.~$\min(w) =x^{-1} \psi_x (\min(w'))$ if $x \ne a$). For short, let $m$, $m'$ denote the respective words $\min(w)$, $\min(w')$. The induction step will consist in constructing some word $u'$ such that inequality \eqref{e2} holds for $w'$, $u'$.

For any acceptable pair $\pi =(a,<)$ with $a \in \Alph (w)$, let $h= h(\pi)$ be maximal such that $au_h$ is a prefix of $m$, and let $H$ be the largest $h(\pi)$ for all such pairs $\pi$. As $u_H \in F(w)$ and begins with $x$, we have $u_H=\psi_x (v)$ for some word $v$. 

Now consider an acceptable pair $\pi =(a,<)$ as above with $h<H$. If $au_h=m$ then we see that $av_q =m'$ for some $q$. Otherwise there exist letters, $y < z$ such that $au_{h+1}=au_h y$ and $ m_{h+2} =au_h z$; whence easily $av_{q+1}=av_qy$ and $m'_{q+2}= av_qz$, and thus $av_{|m'|-1} <m'$.
Now, for any pair $(a,<)$ such that $h=H$ we have either $au_H =m$ or $au_{H+1}=au_Hy < m_{H+2} = m_{H+1}z$, for some letters $y<z$; whence $av=m'$ or $avy < m$. 

Consequently we can take either $u'=v$ or  $u'=vy$. This is the induction step. Clearly $|w'|=\ell' <\ell= |w|$ unless $|\Alph(w)|=1$, a trivial case.

\noindent\emph{Case $2$}: $w$ does not begin with $x$. In this case, we have $w = x^{-1}\psi_x(w')$ for some word $w'$ that does not begin with $x$. Consider $W=xw = \psi_x(w')$. Then, for any acceptable pair $(a,<)$ with $a \neq x$, we have easily $\min(W) = \min (w)$. The same holds if $a=x$ and $aa$ occurs in $w$ because in this case $\min(W)$ begins with $aa$ and $W$ begins with $ay$ for some $y \neq x$; thus $\min(W) \in F(w)$. If $x=a$ and $xx \not \in F(w)$, then the letter $x$ does not occur in $w'$, so inequality \eqref{e2} is trivially satisfied for $w'$ (as stated previously). Thus we can use $W=xw$ instead of $w$ for performing the induction step as in Case 1, ignoring acceptable pairs of the form $(x,<)$.  However, as $|W| =|w| +1$, it is possible that $|w'|=|w|$ or $|w'|=|w| +1$, which are trivial cases corresponding to words $w'$ of the form $yx^p$ for some letter $y \ne x$ and $p \in \NN$. 
\qed\end{pf*}

\begin{example} Recall the finite word $w= baabacababac$ from Example \ref{Ex:min}. 
For the different orders on $\{a,b,c\}$, we have
\begin{itemize}
\item $a<b<c$ or $a<c<b$: $\min(w) = aabacababac$,
\item $b<a<c$ or $b<c<a$: $\min(w) = babac$,
\item $c<a<b$ or $c<b<a$: $\min(w) = cababac$. 
\end{itemize}
It can be verified that a finite word $u$ satisfying \eqref{e2} must begin with $aba$ and one possibility is $u=abacaaaaaa$; thus $w$ is a finite episturmian word. 
\end{example}
\noindent\textit{Note.} In the above example, any two acceptable pairs involving the same letter give the same $\min(w)$, which is not the case in general.

A corollary of Theorem \ref{t3} is the following new characterization of finite Sturmian words (i.e., finite balanced words).

\begin{cor}\label{cor2} A finite word $w$ on $\cA=\{a,b\}$, $a<b$, is not Sturmian (in other words, not  balanced)  if and only if there exists a finite word $u$ such that $aua$ is a prefix of $\min(w)$ and $bub$ is a prefix of $\max(w)$. 
\qed
\end{cor}

\begin{example}
For $w = ababaabaabab$, $\min(w)=aabaabab$ and $\max(w)=babaabaabab$. The  longest common prefix of $a^{-1}\min(w)$ and $b^{-1}\max(w)$ is $abaaba$, which is followed by $b$ in $\min(w)$ and $a$ in $\max(w)$. Thus $w$ is Sturmian. However, if we take $w= aabababaabaab$ for instance, then $w$ is not Sturmian since $\min(w)= auab$ and $\max(w)= bubaabaab$ where $u = aba$.
\end{example}

\begin{remq}
An unrelated connection between finite balanced words (i.e., finite Sturmian words) and lexicographic ordering was recently studied by Jenkinson and Zamboni \cite{oJlZ04char}, who presented three new characterizations of `cyclically' balanced finite words via orderings. Their  characterizations are based on the ordering of  a {\em shift orbit}, either lexicographically or with respect to the $1$-\emph{norm}, which counts the number of occurrences of the symbol $1$ in a given finite word over $\{0,1\}$.
\end{remq}

\section{A characterization of infinite episturmian words in a wide sense} \label{S:infinite}

In this last section, we characterize by lexicographic order the set of all infinite words whose factors are (finite) episturmian. Such infinite words are exactly the episturmian and episkew words, as shown in Proposition \ref{P:wide} below.

\begin{defn} An infinite word is said to be \textbf{\em episturmian in the wide sense} if all of its factors are (finite) episturmian.
\end{defn}

We have the following easy result:

\begin{prop} \label{P:wide} An infinite word is episturmian in the wide sense if and only if it is episturmian or episkew.
\end{prop}
\vspace{-15pt}\begin{pf} 
Let $\bt$ be an infinite word. First suppose that $\bt$ is episturmian in the wide sense. Clearly, if $\bt$ is recurrent, then $\bt$ is episturmian ({\it cf.} proof of (i) $\Rightarrow$ (ii) in Proposition \ref{P:skew}). On the other hand, if $\bt$ is non-recurrent, then $\bt$ is episkew, by Proposition \ref{P:skew}.

Conversely, if $\bt$ is episturmian or episkew, then all of its factors are (finite) episturmian, and hence $\bt$ is episturmian in the wide sense.
\qed\end{pf}

\begin{remq}
Recall that in the $2$-letter case the balanced infinite words (all of whose factors are finite Sturmian) are precisely the Sturmian and skew infinite words. As such, `episturmian words in the wide sense' can be viewed as a natural generalization of balanced infinite words to an arbitrary finite alphabet.
\end{remq}

As a consequence of Theorem \ref{t3}, we obtain the following characterization of episturmian words in the wide sense (episturmian and episkew words).

\begin{cor}\label{cor1*} An infinite word $\bt$ on $\cA$ is episturmian in the wide sense if and only if there exists an infinite word $\bu$ such that 
\begin{equation}a\bu \le \min(\bt) \label{e1*}\end{equation}
for any acceptable pair  $(a, <)$. 
\end{cor}

\vspace{-15pt}\begin{pf} IF part: Inequality \eqref{e1*} holds. So, for any factor $w$ of $\bt$ and any acceptable pair $(a,<)$, we have 
\[
  a\bu_{|m|-1} \leq m \quad \mbox{where $m = \min(w)$}.
\] 
Therefore, by Theorem \ref{t3}, $w$ is a finite episturmian word; whence $\bt$ is episturmian in the wide sense since any factor of $\bt$ is (finite) episturmian.

\smallskip

\noindent  ONLY IF part: $\bt$ is episturmian in the wide sense, so all of its factors are (finite) episturmian; in particular, any prefix $\bt_q$ of $\bt$ is finite episturmian. Therefore, by Theorem \ref{t3},  there exists a finite word, say $u(q)$, such that, for any acceptable pair $(a,<)$, we have
\[
  au(q)_{|m(q)|-1} \leq \min(\bt_q) \quad \mbox{where $m(q) = \min(\bt_q)$}.
\]
 
On the other hand, for any $k \in \NN$ there exists a number $r(k) \in \NN$ such that, for any $q \geq r(k)$, $\bt_q$ contains all the $\min(\bt | k)$ as factors for all acceptable pairs $(a, <)$. It follows then that $\min(\bt | k)$ is a prefix of $\min(\bt_q)$; in particular $|\min(\bt_q)| \ge k$, and hence $|u(q)| \geq k-1$. Thus, the $|u(q)|$ are unbounded.
 
Let us denote by $\bu$ a limit point of the $u(q)$. Then, for any $n$, infinitely many $u(q)$ have $\bu_n$ as a prefix.
 
Now, for any given $k\in \NN ^+$ and acceptable pair $(a,<)$, there exists a $q$ (as above) such that
\[
a\bu _{k-1} = au(q)_{k-1}\le \min(\bt_q)_k = \min(\bt | k).
\]
Thus $a \bu \le \min (\bt)$.
\qed\end{pf}

In the 2-letter case, we have the following characterization of balanced infinite words; in other words, all Sturmian and skew infinite words.

\begin{cor}\label{cor1} An infinite word $\bt$ on $\{a,b\}$, $a<b$, is balanced (i.e., Sturmian or skew) if and only if there exists an infinite word $\bu$ such that 
\begin{equation*}a\bu \le \min(\bt) \le \max(\bt) \le b\bu. \label{e1} \qed\end{equation*} 
\end{cor}

\begin{remq} A variation of the above result appears, under a different guise, in a paper by S.~Gan \cite[Lemma 4.4]{sG01stur}.
\end{remq}


 \bigskip\bigskip

\noindent {\bf GLEN, Amy}* [Corresponding Author] \\
LaCIM, Universit\'{e} du Qu\'{e}bec \`{a} Montr\'{e}al, \\
Case postale 8888, succursale Centre-ville, \\ 
Montr\'{e}al (QC) CANADA, H3C 3P8 \bigskip


\noindent {\bf JUSTIN, Jacques} \\
LIAFA, Universit\'e Paris VII, \\
case 7014, 2 place Jussieu, \\ 75251 Paris Cedex 05, FRANCE \bigskip
 
\noindent {\bf PIRILLO, Giuseppe} \\ 
IASI CNR, Unit\`a di Firenze, \\
Viale Morgagni 67/A, \\ 
50134 Firenze, ITALY

\noindent OR

\noindent Universit\'e de Marne-la-Vall\'ee, \\
5 boulevard Descartes Champs sur Marne, \\ 
77454 Marne-la-Vall\'ee Cedex 2, FRANCE


\begin{thebibliography}{99}
\bibliographystyle{elsart-num-sort}
\footnotesize

\bibitem{jAjS03auto} J.-P.~Allouche, J.~Shallit,  \emph{Automatic Sequences: Theory,
               Applications, Generalizations},  \emph{Cambridge University Press}, UK, 2003.
               

\bibitem{pAgR91repr} P.~Arnoux, G.~Rauzy, Repr\'{e}sentation 
g\'{e}om\'{e}trique de suites de complexit\'{e} $2n+1$, 
\emph{Bull. Soc. Math. France} \textbf{119} (1991), 199--215.


\bibitem{jBpS02stur} J.~Berstel, P.~S\'{e}\'{e}bold,
Sturmian words, in: \emph{M.~Lothaire, Algebraic Combinatorics On
Words, Encyclopedia of Mathematics and its Applications}, vol. 90,
\emph{Cambridge University Press}, U.K., 2002, pp.~45--110.


\bibitem{aD97stur} A.~de~Luca, Sturmian words: structure,
combinatorics and their arithmetics, \emph{Theoret. Comput. Sci.}
\textbf{183} (1997), 45--82, doi:10.1016/S0304-3975(96)00310-6.


\bibitem{xDjJgP01epis} X.~Droubay, J.~Justin, G.~Pirillo,
Episturmian words and some constructions of de Luca and Rauzy,
\emph{Theoret. Comput. Sci.} \textbf{255} (2001), 539--553,
doi:10.1016/S0304-3975(99)00320-5.


\bibitem{xDgP99pali} X.~Droubay, G.~Pirillo, Palindromes and 
{Sturmian} words,
\emph{Theoret. Comput. Sci.} \textbf{223} (1999), 73--85,
doi:10.1016/S0304-3975(97)00188-6.


\bibitem{sG01stur} S. Gan, Sturmian sequences and the lexicographic world, 
{\it Proc. Amer. Math. Soc.} {\bf 129} (2001), 1445--1451.


\bibitem{aG05powe} A.~Glen, Powers in a class of $\mathcal{A}$-strict standard 
               episturmian words, in: \emph{$5$th International Conference on Words, Universit\'{e} du Qu\'{e}bec \`{a} Montr\'{e}al}, Publications du LaCIM \textbf{36} (2005), 249--263.  {\em Theoret. Comput. Sci.} (in press), doi:10.1016/j.tcs.2007.03.023.
                   
                   
\bibitem{aG06acha} A.~Glen, A characterization of fine words over a finite alphabet, {International School and Conference on Combinatorics, Automata and Number Theory, Universit\'{e} de Li\'{e}ge, Belgium} (2006), pp. 9.  {\em Theoret. Comput. Sci.} (accepted).                


\bibitem{aHrT00char} A.~Heinis, R.~Tijdeman, Characterisation of asymptotically {S}turmian sequences, \emph{Publ. Math. Debrecen} {\bf 56} (3--4) (2000), 415--430. 


\bibitem{oJlZ04char} O.~Jenkinson, L.Q.~Zamboni, Characterisations of balanced words via orderings, \emph{Theoret. Comput. Sci} \textbf{310} (2004), 247--271, doi:10.1016/S0304-3975(03)00397-9. 


\bibitem{jJgP97deci} J.~Justin, G.~Pirillo, Decimations and {S}turmian words, \emph{Theor. Inform. Appl.} {\bf 31} (3) (1997), 271--290.
	
	
\bibitem{jJgP02epis} J.~Justin, G.~Pirillo, Episturmian words and episturmian morphisms, \emph{Theoret. Comput. Sci.} \textbf{276} (2002), 281--313, doi:10.1016/S0304-3975(01)00207-9.


\bibitem{jJgP02onac} J.~Justin, G.~Pirillo,
   On a characteristic property of {Arnoux-Rauzy} sequences, \emph{Theor. Inform. Appl.} \textbf{36} (4) (2002),  385--388.  


\bibitem{jJgP04epis} J.~Justin,  G.~Pirillo,  Episturmian words: shifts, morphisms and numeration systems, \emph{Internat. J. Found. Comput. Sci.} \textbf{15} (2) (2004), 329--348, doi:10.1142/S0129054104002455. 


\bibitem{jJlV00retu} J.~Justin, L.~Vuillon,
Return words in {Sturmian} and episturmian words, \emph{Theor. Inform. Appl.} \textbf{34} (5) (2000), 343--356.

\bibitem{fMlZ02onth} F.~Mignosi, L.Q.~Zamboni, On the number of {A}rnoux-{R}auzy words, \emph{Acta Arith.} \textbf{101} (2) (2002), 121--129.


\bibitem{gHmM40symb} M.~Morse, G.A.~Hedlund, Symbolic dynamics
II. Sturmian trajectories, \emph{Amer. J. Math.} \textbf{62}
(1940), 1--42.


\bibitem{gP99anew} G.~Pirillo, A new characteristic property of the palindrome prefixes of a standard Sturmian word,  {\it S\'em. Lothar. Combin.} {\bf 43} (1999), pp. 3.


\bibitem{gP05ineq} G.~Pirillo, Inequalities characterizing standard {S}turmian and
              episturmian words, \emph{Theoret. Comput. Sci.} \textbf{341} (1--3) (2005),  276--292, doi:10.1016/j.tcs.2005.04.008.
              
              
\bibitem{gP05mors} G.~Pirillo,  Morse and {H}edlund's skew {S}turmian words revisited,  \emph{Ann. Comb.} (accepted). 


\bibitem{nP02subs} N. Pytheas Fogg, Substitutions in Dynamics, Arithmetics and Combinatorics, Lecture Notes in Mathematics, vol. 1794, \emph{Springer-Verlag}, Berlin, 2002. 

              

\bibitem{gR03conj} G.~Richomme, Conjugacy and episturmian morphisms, \emph{Theoret. Comput. Sci.} {\bf 302} (2003), 1--34, doi:10.1016/S0304-3975(02)00726-0.
	
   
\bibitem{rRlZ00agen} R.N.~Risley, L.Q.~Zamboni, 
A generalization of {Sturmian} sequences: {Combinatorial}
structure and transcendence, \emph{Acta Arith.} \textbf{95} (2)
(2000), 167--184.


\end{thebibliography}
\end{document}